\UseRawInputEncoding
\documentclass[12pt, leqno]{article}
\usepackage{latexsym}
\usepackage{enumerate,amsmath,amsthm,amssymb,amsfonts,amscd}

\newtheorem{defn}{Definition}[section]
\newtheorem{thm}[defn]{Theorem}
\newtheorem{cor}[defn]{Corollary}
\newtheorem{prop}[defn]{Proposition}
\newtheorem{lem}[defn]{Lemma}

\numberwithin{equation}{section}

\def\F{\mbox{${\cal F}$}}

\newcommand{\cc}{\mathbb C}
\newcommand{\zz}{\mathbb Z}
\newcommand{\nn}{\mathbb N}

\input{arrow.tex}

\begin{document}
\title{\bf On an affirmative solution to the non-commutative Michael's problem, with applications to automatic continuity theory}
\author{S. R. PATEL}
\maketitle \vspace{1.5cm}

\noindent Department of Mathematics, C. U. Shah University, Wadhwan City,\\
Gujarat, INDIA\\
E-mail: srpatel.math@gmail.com
\begin{tabbing}
\> Secondary: 13F25, 32A10, 46A32,\\
\> 46E25, 46H10, 46H40
\end{tabbing}
\newpage
\noindent {\bf Abstract.} In 1952, Michael posed a question about
functional continuity of commutative Fr\'{e}chet algebras in his
memoir, known as Michael's problem in the literature. We settle
this in the affirmative for the non-commutative Fr\'{e}chet
algebras, along with its equivalent form. Indeed, we continue our
previous works to directly attack these problems. We show that the
existence of a discontinuous character on the non-commutative test
case for this problem -- the Fr\'{e}chet algebra ${\cal U}_{nc}$
of all entire functions on $\ell^\infty$ -- allows us to generate
another Fr\'{e}chet algebra topology on ${\cal U}_{nc}$,
inequivalent to the usual Fr\'{e}chet algebra topology. An
important tool is a topological version of the tensor algebra over
a Banach space. The elementary, but crucial, idea is to express
the test algebra ${\cal U}_{nc}$ as a {\it weighted Fr\'{e}chet
tensor algebra} over the Banach space $\ell^1$.

We give some important applications in automatic continuity
theory. In particular, we provide affirmative answers to
significant, long-standing problems from 1967 (extending Johnson's
result for non-commutative, semi-simple Fr\'{e}chet algebras;
Dales' survey, Question 9), 1976 (extending Esterle's result,
giving another approach for Michael's problem for more general
$(F)$-algebras), and 1978 (a weaker version of an equivalent form
of the non-commutative Michael's problem; Dales' survey, p. 143).
\newpage
\section{Introduction and statement of the main theorem}
{\bf History of Michael's problem.}

An important subject of study in the theory of Fr\'{e}chet
algebras is certain questions from automatic continuity theory.
Viz., we may consider the (non-)uniqueness of the Fr\'{e}chet
algebra topology (briefly: the Fr\'{e}chet topology in the sequel)
on certain Fr\'{e}chet algebras, (dis)continuity property of
certain homomorphisms between certain Fr\'{e}chet algebras as well
as that of (higher point) derivations on Fr\'{e}chet algebras.
This theory has applications to non-commutative rings and algebras
(e.g., the non-commutative Singer-Wermer conjecture in Fr\'{e}chet
algebras ~\cite{33, 37, 5, 31}).

In 1952, Michael posed a question about functional continuity of
commutative Fr\'{e}chet algebras ~\cite{27}. This problem has a
very strong historical background. A number of significant
analysts worked on this problem since 1952. The non-commutative
analogues of Michael's problems also exist (see below). Some
analysts discussed various test cases and approaches (e.g., Dixon
and Esterle ~\cite{15}, Dales, author and Read ~\cite{14}). Others
discussed various equivalent forms (e.g., Craw~\cite{9}, Dixon and
Fremlin ~\cite{16}, Akkar~\cite{1}).

Despite a lot of efforts by analysts to solve Michael's problem,
it seems that only six significant ideas appeared in the
literature since 1952 ~\cite{3, 34, 8, 15, 18, 14}. The strong
partial result was obtained by Arens in 1958. As far as we know,
the most latest effort was made by Dales, author and Read in 2010.
Dixon and Esterle discussed the test case ${\cal U}_{nc}$ (the
non-commutative analogue of $\cal U$) for the non-commutative
Michael's problem in ~\cite{15}.

\noindent {\bf Our past work as a background.}

So far, we have established the uniqueness of the Fr\'{e}chet
topology for {\it all} Fr\'{e}chet algebras of power series
(briefly: FrAPS) in one indeterminate ~\cite{28} and in finitely
many indeterminates ~\cite{29}. Dales, author and Read have also
obtained an affirmative solution (in a stronger form) to
Dales-McClure's problem (1977) among other interesting automatic
continuity results, including a reduction of Michael's problem
~\cite[Theorems 11.2 and 10.1 (ii)]{14}. We have used
discontinuity of derivations to give other inequivalent
Fr\'{e}chet topology to certain Fr\'{e}chet algebras, solving
Loy's problem (1974) in the Fr\'{e}chet case ~\cite{30}. Also, by
using the Read's method, we have constructed two maiden examples
of Fr\'{e}chet algebras ($\F_\infty \,\textrm{and}\,
\F_\infty\,\oplus\,\F_\infty$) with countably many mutually
inequivalent Fr\'{e}chet topologies, to discuss various avatars of
the Singer-Wermer conjecture in Fr\'{e}chet algebras ~\cite{33,
31}. Such examples are not known even in the Banach algebra case.
As we shall see below, {\it all} these works have turned out to be
stepping stones for the complete solution to the non-commutative
Michael's problem (e.g., see ~\cite[Question, p. 135]{28} and
~\cite[Theorem 11.2]{14}).

Since 1958, it was not clear to functional analysts how to extend
the Arens' result from the finitely generated case to the
countably generated case, in order to obtain the functional
continuity of the test algebra $\cal U$ (resp., of ${\cal
U}_{nc}$). One strongly believed that there must have some ways to
apply appropriate methods within automatic continuity theory as
the problem falls in this theory. So, we, here, represent the
second approach for the non-commutative analogue of Michael's
problem; in another paper, we exclusively discuss the first
approach for the commutative case. In both the approaches, an
important tool is a topological version of the tensor (symmetric)
algebra over a Banach space. The elementary, but crucial, idea is
to express the test algebra(s) ${\cal U}_{nc}$ (and $\cal U$) as a
{\it weighted Fr\'{e}chet tensor (symmetric) algebra} over the
Banach space $\ell^1$ ~\cite{14, 30}.

\noindent {\bf Preliminaries.}

Throughout the paper, ``algebra" will mean a non-zero, complex
(non-)commutative algebra with identity. A topological algebra $A$
is functionally continuous if every character on $A$ is
continuous, and $A$ is functionally bounded if every character on
$A$ is bounded on bounded sets. A {\it Fr\'{e}chet algebra} is a
complete, metrizable locally convex algebra $A$ whose topology
$\tau$ may be defined by an increasing sequence $(p_m)_{m\in \nn}$
of submultiplicative seminorms. The basic theory of Fr\'{e}chet
algebras was introduced in ~\cite{27, 19}. If $A$ is a Fr\'{e}chet
algebra, then it is functionally continuous if and only if it is
functionally bounded ~\cite[Theorem 1]{16}. However, there are
commutative, complete LMC-algebras which are functionally bounded,
but not functionally continuous; e.g., see ~\cite[pp. 52-53]{27}.

We write $\cal B$ for the algebra $\cc_{nc}[[X_0, X_1, \dots]]$ of
all formal power series in countably many non-commuting
indeterminates $X_0, X_1, \dots$, with complex coefficients. We
shall also require the commutative analogue of $\cal B$, denoted
by $\F_\infty$ = $\cc[[X_0, X_1, \dots]]$. A fuller description of
these algebras is given in ~\cite[\S 9]{14}; see also ~\cite{18}
for $\F_\infty$. Both are Fr\'{e}chet algebras under the usual
topology $\tau_c$ of coordinatewise convergence defined by the
coordinate projections $\pi_I, I \in S_{nc}$ the free semigroup in
countably many non-commuting indeterminates $X_0, X_1, \dots$ for
$\cal B$ (resp., $\pi_I, I \in (\zz^+)^{<\omega}$ for
$\F_\infty$). The defining sequence $(p_m^{'})$ of seminorms is:
$$p_m^{'}(\sum_{J \in S_{nc}}\lambda_J X^{\otimes J}) = \sum\{\mid
\lambda_J \mid : J \in (\zz^+)^m, \textrm{rank}\,J \leq m\}\;(m
\in \nn)$$ (resp., $p_m^{'}(\sum_{J \in
(\zz^+)^{<\omega}}\lambda_J X^J) = \sum\{\mid \lambda_J \mid \,: J
\in (\zz^+)^m, \mid J\mid \,\leq m\}\;(m \in \nn)$). Clearly, each
$p_m^{'}$ is a proper seminorm on $\cal B$ (resp., on
$\F_\infty$).

The algebra $\cal B$ is a graded algebra. We shall also require in
a future proof the 'averaging map' and symmetrizing map
$\tilde{\sigma}$ on $\cal B$, and symmetric elements of $\cal B$.
There is a product in ${\cal B}_{\textrm{sym}}$, and with this
product $\vee$, it is a commutative, unital algebra and ${\cal
B}_{\textrm{sym}} = \epsilon(\F_\infty)$, where $\epsilon:
\F_\infty \,\rightarrow \cal B$ is a continuous linear embedding
~\cite{14}. In fact, $\epsilon$ is a continuous, injective
homomorphism as we shall see in \S 3.

\noindent {\bf Main theorem, with comments on the strategies.}

In this paper, we shall be concerned with the affirmative
solutions to Michael's two problems. In \S 2, we shall briefly
discuss some general theory of topological tensor (symmetric)
algebras over a Banach space, in order to establish notation that
will be used throughout the paper. In particular, we mainly study
certain semigroup (Banach/Fr\'{e}chet) algebras over the
semigroups $S = (\zz^+)^{<\omega}$ and $S_{nc}$, which are graded
subalgebras of $\F_\infty$ and $\cal B$ ~\cite[\S 9, \S 10]{14}.
These algebras are required for our approach to solve Michael's
problem in \S 3. Thus these algebras over the Banach space
$\ell^1$ are the technical main-spring of the paper.

In our approach, the essential ingredient is as follows. First, we
describe the test algebra ${\cal U}_{nc}$ in terms of weighted
Fr\'{e}chet tensor algebra $\widehat{\bigotimes}_WE$ over the
Banach space $E\,=\,\ell^1$ (the commutative analogue of this
example ~\cite{30}, Question 3 of ~\cite{28} and ~\cite{31, 33},
motivated the ideas of the present paper). Then we show that the
existence of a discontinuous character $\phi$ on ${\cal U}_{nc}$
allows us to generate another Fr\'{e}chet topology $\tau$ on
${\cal U}_{nc}$ (resp., on $\cal U$), inequivalent to the usual
Fr\'{e}chet topology $\tau_0$ on these algebras. This contradicts
the fact that ${\cal U}_{nc}$ (resp., $\cal U$), being a
non-(commutative), semi-simple FrAPS in $\cal B$ (resp., in
$\F_\infty$), has a unique Fr\'{e}chet topology (see Corollary 4.3
below; resp., ~\cite{6, 28, 29, 14}). We give some crucial remarks
on our approach. We have the following

\noindent {\bf Main Theorem.} {\it The algebra $({\cal U}_{nc},
\tau)$ is complete with respect to $(\|\cdot\|_n)$. The derivation
$\partial/\partial\,X_0\,:({\cal
U}_{nc},\,\tau)\,\rightarrow\,({\cal U}_{nc},\,\tau)$ is
discontinuous, and its separating subspace is all of ${\cal
U}_{nc}$. The derivation $\partial/\partial\,X_0\,:\,({\cal
U},\,\tau)\,\rightarrow\,({\cal U},\,\tau)$ is also discontinuous,
and its separating subspace is all of ${\cal U}$.} $\hfill \Box$

We obtain the functional continuity of the test algebra $({\cal
U}_{nc},\,\tau_0)$ in the following

\noindent {\bf Corollary.} {\it All characters on the Fr\'{e}chet
algebra $({\cal U}_{nc}, \tau_0)$ are continuous. In particular,
these characters are, indeed, the point evaluation mappings, that
is, if $\phi$ is a character on ${\cal U}_{nc}$, then there exists
$z\,\in\,\ell^\infty$ such that $\phi\,=\phi_z$, where
$\phi_z(f)\,=\,f(z)$ for all $f\,\in\,{\cal U}_{nc}$.} $\hfill
\Box$

As consequences, we have affirmative solutions to Michael's two
problems as follows (see Corollary 3.4 and Corollary 3.5 below).

\noindent {\bf Statements.} {\it 1. Every non-commutative
Fr\'{e}chet algebra is functionally continuous. 2. Every
non-commutative, complete LMC-algebra is functionally bounded.}

In \S 3 and \S 4, we summarize the previous developments on
Michael's problems in view of our affirmative solution. We shall
discuss the non-commutative analogue of the Carpenter's result
~\cite[Corollary 6]{6}. It is a known fact that another equivalent
form of Michael's problem is the long-standing problem of
continuity of a homomorphism from a Fr\'{e}chet algebra $B$ into a
commutative semi-simple Fr\'{e}chet algebra $A$ ~\cite[p.
143]{10}. We shall discuss a weaker version of this equivalent
form in the non-commutative case. In fact, it is a surprising
consequence that attempts to solve the non-commutative analogues
of Michael's problems, lead us to certain important automatic
continuity results including the uniqueness of the Fr\'{e}chet
topology for non-commutative, semi-simple Fr\'{e}chet algebras
(extending Johnson's result (1967)); see Corollary 4.3 below. This
also leads us to another approach to answer Michael's problem for
more general complete, metrizable topological algebras (briefly:
$(F)$-algebras) in both the cases (commutative as well as
non-commutative) by extending Esterle's result ~\cite{17}; many
thanks to Thomas stability lemma ~\cite{35}.
\section{Topological tensor (symmetric) algebras over a Banach space} In \S 3, our important tool is a topological version of the tensor
(symmetric) algebra over a Banach space appearing in ~\cite{23}.
However, a fuller description of what we require is given in
~\cite{13, 11}. We briefly recall some notation, which will be
used throughout the paper. Let $E = \hat{\bigotimes}^1E$ be a
Banach space, and take $\hat{\bigotimes}^0E = \cc$. For each $p
\geq 2$, write $\hat{\bigotimes}^pE$ (resp.,
$\check{\bigotimes}^pE$) for the completion of $\bigotimes^pE$
with respect to (briefly: w.r.t.) the projective tensor product
norm $\|\cdot\|_\pi$ (resp., the equicontinuous tensor product
norm $\|\cdot\|_\epsilon$). When it is unnecessary to distinguish,
we may write $\hat{\bigotimes}^pE$ for either completion, and
$\|\cdot\|$ for the specified norm. Then $\widehat{\bigotimes}E$
is a non-commutative, unital Fr\'{e}chet tensor algebra over $E$
(if $\textrm{dim}\,E > 1$) w.r.t. the product
\begin{equation} \label{tensor}
(\sum_pu_p)\,\otimes\,(\sum_pv_p)\,=\,\sum_p(\sum_{i+j=p}u_i\,\otimes\,v_j)
\end{equation}
and the coefficientwise convergence topology defined by an
increasing sequence $(\|\cdot\|_m)$ of seminorms, where
$\|\sum_p\,u_p\|_m\,=\,\sum_{p=0}^m\|u_p\|$. We refer to the
subspaces $\hat{\bigotimes}^pE$ as the homogeneous subspaces of
$\widehat{\bigotimes}E$.

We shall also need in a future proof a commutative analogue
$\widehat{\bigvee}E$ of $\widehat{\bigotimes}E$, and a closed,
linear subspace $\hat{\bigvee}^pE$ of $\hat{\bigotimes}^pE$. Then
$\widehat{\bigvee}E$ is a commutative, unital Fr\'{e}chet
symmetric algebra over $E$ w.r.t. the product
\begin{equation} \label{symmetric}(\sum_pu_p)\,\vee\,(\sum_pv_p)\,=\,\sum_p(\sum_{i+j=p}u_i\,\vee\,v_j)
\end{equation}
and the same coefficientwise convergence topology defined by
$(\|\cdot\|_m)$. We remark that if $\tilde{\sigma} \,=\,\bigoplus
\tilde{\sigma}_{p}$ is the continuous, symmetrizing epimorphism
$\tilde{\sigma} : \widehat{\bigotimes}E \rightarrow
\widehat{\bigvee}E$, then $\textrm{ker}\,\tilde{\sigma}$ is the
closed, two-sided ideal of $\widehat{\bigotimes}E$ generated by
$$\{u\otimes v - v\otimes u : u,\, v \in \hat{\bigotimes}^1E =
E\},$$ so that the algebras in ~\cite{23} and the algebras
$\widehat{\bigvee}E$ are topologically isomorphic ~\cite[p.
315]{13}.

There are some important Banach subalgebras of the above examples
~\cite[Example 2.2.46 (ii)]{11}. Let $\omega$ be a weight on
$\zz^+$ and let $A$ be one of the algebras $\widehat{\bigotimes}E$
and $\widehat{\bigvee}E$. (When it can cause no confusion, the
symbol $(\zz^+)^{<\omega}$ is used to denote the subsemigroup
consisting of all $\zz^+-$valued sequences that are eventually
$0$.) Define $\widehat{\bigotimes}_{\omega}E$ and
$\widehat{\bigvee}_{\omega}E$, respectively, as
\begin{equation} \label{Banach}\{u\,=\,\sum_{p=0}^\infty\,u_p\,\in\,A\,:\,\|u\|_\omega\,:=\,\sum_{p=0}^\infty\|u_p\|\omega(p)\,<\,\infty\}.
\end{equation}
We obtain two unital Banach subalgebras, with
$\widehat{\bigvee}_{\omega}E$ being commutative. These algebras
are called weighted Banach tensor algebras and weighted Banach
symmetric algebras, respectively. The algebra
$\widehat{\bigotimes}_{\omega}E$, where $\omega \equiv 1$ on
$\zz^+$, is the Banach tensor algebra which is the starting point
for the constructions in ~\cite{23}. We recall the following
theorem from ~\cite{13} (also, cf. ~\cite[Satz 2 and Satz 3]{23}),
describing the maximal ideal space and semi-simplicity of
$\widehat{\bigvee}_{\omega}E$.
\begin{thm} \label{Thm 1.2_DM} \begin{enumerate}
\item[{\rm (i)}] The space of characters on
$\widehat{\bigvee}_{\omega}E$ is homeomorphic with the ball
$\{\lambda\,\in\,E{'}\,:\,\|\lambda\|\,\leq\,\inf_{p}\omega(p)^\frac{1}{p}\,\equiv\,\omega_\infty\}\;\;(w^*-\textrm{topology}).$
\item[{\rm (ii)}] $\widehat{\bigvee}_{\omega}E$ (w.r.t. the
equicontinuous tensor norm $\|\cdot\|_\epsilon$) is semi-simple if
and only if $\omega_\infty\,>\,0$. \item[{\rm (iii)}]
$\widehat{\bigvee}_{\omega}E$ (w.r.t. the projective tensor norm
$\|\cdot\|_\pi$) is semi-simple if and only if
$\omega_\infty\,>\,0$ and $E$ has the approximation property.
$\hfill \Box$
\end{enumerate}
\end{thm}
Along the same lines, for an increasing sequence
$W\,=\,(\omega_m)$ of weights on $\zz^+$, we define
$\widehat{\bigotimes}_{W}E$ and $\widehat{\bigvee}_{W}E$,
respectively, as
\begin{equation} \label{Frechet} \{u\,=\,\sum_{p=0}^\infty\,u_p\,\in\,A\,:\,p_m(u)\,:=\,\sum_{p=0}^\infty\|u_p\|\omega_m(p)\,<\,\infty\;\textrm{for
all}\, m\,\in\,\nn\}. \end{equation} We obtain two unital
Fr\'{e}chet subalgebras, with $\widehat{\bigvee}_{W}E$ being
commutative. These algebras are called weighted Fr\'{e}chet tensor
algebras and weighted Fr\'{e}chet symmetric algebras,
respectively. Clearly, from \eqref{Frechet}, their Arens-Michael
representations are given by $\widehat{\bigvee}_{W}E = \bigcap_{m
\in \nn}(\widehat{\bigvee}_{\omega_m}E,\, p_m)$ and
$\widehat{\bigotimes}_{W}E$
$=\,\bigcap_{m\,\in\,\nn}(\widehat{\bigotimes}_{\omega_m}E,\,
p_m)$, where $(p_m)$ is a sequence of norms on the respective
Fr\'{e}chet algebra.

Set $E=\ell^1(\zz^+)$, the Banach space for the remaining part of
this section. Then, as in ~\cite[\S 10]{14}, for each $p \in \nn$,
$\hat{\bigotimes}^pE$ can be identified with $\ell^1((\zz^+)^p)$
as a Banach space. This Banach space can also be viewed as the
space of absolutely summable functions on $(\zz^+)^p$
~\cite[Chapter 45]{38} (equivalently, the space $A^+(D^p)$ of
functions in $A(D^p)$ ($D^p$ the product of $p$-copies of $D$ in
$\cc$), with absolutely convergent Taylor series on $\Gamma^p$
(the product of $p$-copies of the unit circle $\Gamma$ in $\cc$),
w.r.t. the $\|\cdot\|_{\ell^1}$-norm on the space).
\begin{thm} \label{p 315_DM2} For each $p \in \nn$, the Banach space $\hat{\bigotimes}^pE$ is a Banach algebra, and the character
space of $\hat{\bigotimes}^pE$ is homeomorphic with the polydisc
$D^p$.
\end{thm}
\noindent{\bf Proof.} By ~\cite[Proposition 1.3.11]{11},
$\bigotimes^pE$ is a commutative, unital algebra w.r.t. the
product
\begin{equation}
\label{product-D}(f_1 \otimes \dots \otimes f_p)\cdot(g_1
 \otimes \dots \otimes g_p) = (f_1g_1 \otimes
\dots \otimes f_pg_p)\;(f_i, g_i \in E,\;i \in \nn_p).
\end{equation}
Then, it is easy to see that it is a normed algebra w.r.t. the
projective tensor norm $\|\cdot\|_\pi$ (note that $\|\cdot\|_\pi$
is submultiplicative on $\bigotimes^pE$). Indeed,
$\hat{\bigotimes}^pE$ is a Banach algebra by noticing that the
product can be extended to $\hat{\bigotimes}^pE$ (cf. ~\cite[p.
315]{13}). Thus, for each $p \in \nn$, $\hat{\bigotimes}^pE$ can
be identified with $\ell^1((\zz^+)^p)$ as a Banach algebra via the
homomorphism $$\theta\,:\,(f_1 \otimes \dots \otimes
f_p)(r)\,\mapsto\,f_1(r_1)\cdots f_p(r_p),\; \textrm{where}\;
r\,=\,(r_1,\,\dots,\,r_p)\,\in\,(\zz^+)^p.$$ So, the character
space of $\hat{\bigotimes}^pE$ is homeomorphic with the polydisc
$D^p$ for each $p \in \nn$. $\hfill \Box$

The main point of this identification should be emphasized. The
completion of the tensor power of finite copies of the Banach
space (which is also a singly generated Banach algebra) w.r.t.
$\|\cdot\|_\pi$ is identified with a finitely generated Banach
algebra. Unfortunately, we cannot extend this to $\nn$-fold tensor
product within the framework of Banach algebras. However, we can,
indeed, extend this to $\nn$-fold tensor product within the
framework of Fr\'{e}chet algebras (e.g., the algebra
$\widehat{\bigotimes}E$ and its commutative analogue
$\widehat{\bigvee}E$ below).

As above, $\hat{\bigvee}^pE$ is a closed, linear subspace of
$\hat{\bigotimes}^pE$, consisting of the symmetric elements.
Recall that $\hat{\bigvee}^pE$ is the range of the symmetrizing
map $\tilde{\sigma}_{p}$ (equivalently, the projection of norm
$1$) on $\hat{\bigotimes}^pE$. In fact, every element of
$\hat{\bigotimes}^pE$ is symmetric due to the fact that
$\ell^1((\zz^+)^p)$ is identified with $A^+(D^p)$ and so,
$\partial_kf(z)\;(z\,\in\,D^p)$ will be invariant under
permutations of $1,\,2,\,\dots,\,k$ (cf. ~\cite[p. 320]{13}). So,
$\hat{\bigvee}^pE\,=\,\hat{\bigotimes}^pE$ by ~\cite[Chapter
45]{38} and ~\cite[p. 321]{13}. Thus, as above, for each $p \in
\nn$, $\hat{\bigvee}^pE\,\cong\,\ell^1((\zz^+)^p)$ as a Banach
algebra.

As above, the algebra $\widehat{\bigotimes}E$ is naturally
identified with a graded subalgebra of $\cal B$, and
$\widehat{\bigvee}E$ is naturally identified with a graded
subalgebra of $({\cal B}_{\textrm{sym}}, \vee)$ ~\cite[\S 10]{14},
called the unital Fr\'{e}chet symmetric algebra over $E$. The
latter algebra can also be viewed as the algebra of $p$-absolutely
summable symmetric functions on $(\zz^+)^{<\omega} =
\bigcup_{p\in\nn}(\zz^+)^p$ for all $p \in \nn$. That is, if $f =
\sum_pf_p$ is a symmetric function on $(\zz^+)^{<\omega}$, then
$f_p$ is an absolutely summable function on $(\zz^+)^p$ for each
$p \in \nn$. The former algebra $\widehat{\bigotimes}E$ can also
be viewed as the algebra of $p$-absolutely summable functions on
$(\zz^+)^{<\omega}$.\\

\noindent {\bf Remarks A.} 1. We shall often write the Banach
algebra $\widehat{\bigvee}_{\omega}E$, where $\omega \equiv 1$ on
$\zz^+$, as ${\cal U}_{1} = \ell^1((\zz^+)^{<\omega})$ throughout
the paper (see \S 3 below), since it is isometrically isomorphic
with $\ell^1((\zz^+)^{<\omega})$ ~\cite[p. 142]{14}. So, by
Theorem 2.1, it may be viewed as the algebra of functions in
$A(\ell_{[1]}^\infty)$ with absolutely convergent Taylor series
(equivalently, absolutely summable symmetric functions on
$(\zz^+)^{<\omega}$), where $\ell_{[1]}^\infty$ is the unit
polydisc in $\cc^\infty$; compare ~\cite{8} and the above
situation in the finite-dimensional case.

\noindent 2. Recall that Dales, author and Read gave a very long
proof of showing $\ell^1((\zz^+)^{<\omega})$ a Banach algebra of
power series (briefly: BAPS) ~\cite[Theorem 10.1 (i)]{14}. We
drastically reduce this proof in the following
\begin{lem} \label{Thm 10.1 (i)_DPR} There is a topological isomorphism
$\theta$ of $\ell^1((\zz^+)^{<\omega})$ onto $\ell^1(\zz^+)$, and
so the Banach algebra $\ell^1((\zz^+)^{<\omega})$ is
(isometrically isomorphic to) a semi-simple Banach algebra of
power series $\ell^1(\zz^+)$.
\end{lem}
\noindent {\bf Proof.} First, we note that the homomorphism,
discussed in the proof of ~\cite[Theorem 1.2 (iii)]{13}, from
$\widehat{\bigvee}_{\omega}E$, where $\omega \equiv 1$ on $\zz^+$,
into $(\F, \,\omega)$ (which is, in fact, $\ell^1(\zz^+)$), is,
indeed, an injective epimorphism. Moreover, by ~\cite[\S 9]{22},
it is continuous since $\ell^1(\zz^+)$ is a BAPS in $\F$. So, by
the open mapping theorem, this isomorphism is a topological
isomorphism between ${\cal U}_{1}$ and $\ell^1(\zz^+)$. Thus the
norm $q_1$ of ${\cal U}_{1}$ is equivalent to $\|\cdot\|_{\ell^1}$
on a semi-simple BAPS $\ell^1(\zz^+)$. $\hfill \Box$

However, the difference between the above proof and the proof of
~\cite[Theorem 10.1 (i)]{14} is as follows. In the above proof,
the bounded linear functional $\lambda_k$ on
$\textrm{SE}_\pi(k),\;k\,\in\,\nn$, need not be higher point
derivations of infinite order on $\ell^1(\zz^+)$ (cf. ~\cite[\S 1,
p. 169]{12}, introduction in ~\cite[\S 2]{13} and ~\cite[Theorem
10.1]{14}).
\section{Our approach}
\noindent {\bf Background.}

We now discuss the second approach to solve this problem for
non-commutative Fr\'{e}chet algebras. The method used in the first
approach works only for commutative Fr\'{e}chet algebras, so it is
the need of the hour to develop another approach, which would work
for the non-commutative case as well. However, we show that
solving problem for the non-commutative case would suffice to
solve the problem for the commutative case below. We remark that
if $\phi$ is a discontinuous character on some non-commutative,
unital Fr\'{e}chet algebra $A$, then there is a discontinuous
character on the test case ${\cal U}_{nc}$ ~\cite[Proposition
2.1]{15}; the same statement holds for the test case $\cal U$ as
well ~\cite[REMARK 2.2]{15} and ~\cite[Theorem 2.7]{18}.
\begin{prop} \label{Michael_Esterle} There exists a discontinuous
character on some commutative, unital Fr\'{e}chet algebra if and
only if there exists a discontinuous character on some
non-commutative, unital Fr\'{e}chet algebra.
\end{prop}
\noindent{\bf Proof.} First, we remark that it is sufficient to
prove this proposition for the test cases $({\cal U}, \tau_0)$ and
$({\cal U}_{nc}, \tau_0)$ by ~\cite[Theorem 2.7]{18} (see, also,
~\cite[Theorem 3.4]{10} or ~\cite[Theorem 9]{8} in addition to the
comments above). Let $\phi$ be a discontinuous character on ${\cal
U}_{nc}$. Then $\phi$ does not belong to $M({\cal U}_{nc}) \cong
\ell^\infty$ ~\cite{15, 18}. Since $$M({\cal U}_{nc}) \cong
\ell^\infty \cong M({\cal U})$$ by ~\cite[REMARK 2.2]{15}, there
is a discontinuous character, say $\phi$ again, on $\cal U$. The
same holds true in the reverse direction as well. In fact, we can
associate these discontinuous characters as follows. Recall that
${\cal U}_{nc} \cong \widehat{\bigotimes}_WE$, a weighted
Fr\'{e}chet tensor subalgebra of $\widehat{\bigotimes}E$, is a
graded subalgebra of $\cal B$. Hence, we have a continuous
symmetrizing epimorphism $$\tilde{\sigma} =
\oplus\tilde{\sigma_{p}} : {\cal U}_{nc} \rightarrow {\cal
U}_{nc},\; \textrm{where}\; \tilde{\sigma_{p}} : {\cal
U}_{nc}^{(p)} \rightarrow {\cal U}_{nc}^{(p)}$$ is the averaging
map on ${\cal U}_{nc}^{(p)}$ ~\cite{13, 14}. Clearly, if $\phi$ is
a discontinuous character on $\cal U$, then $\phi \circ
\tilde{\sigma} = \psi$ is a discontinuous character on ${\cal
U}_{nc}$. For the reverse direction, since a continuous
monomorphism $\epsilon|_{{\cal U}} : {\cal U} \rightarrow {\cal
U}_{nc}$ is a natural inverse of $\pi|_{{\cal U} =
\tilde{\sigma}({\cal U}_{nc})}$ (defined in the proof of the main
theorem below), if $\psi$ is a discontinuous character on ${\cal
U}_{nc}$, then $\phi = \psi \circ \epsilon|_{{\cal U}}$ is a
discontinuous character on $\cal U$. $\hfill \Box$

We now discuss our approach, in order to solve Michael's problem
affirmatively. We first see that the test case ${\cal U}_{nc}$
(resp., $\cal U$) is a Fr\'{e}chet algebra of the kind
$\widehat{\bigotimes}_{W}E$ (resp., $\widehat{\bigvee}_{W}E$),
where $E = \ell^1(\zz^+)$ and $W = (\omega_m)$ is an increasing
sequence of weights on $\zz^+$ (note that the germ of the ideas
for our approach lies in ~\cite{30, 28, 31, 33}). Following
~\cite[Definition 9.2 and \S 10]{14}, we have Banach algebras
$\ell^1(S,\, \omega) = {\cal U}_1$, where $\omega \equiv 1$ on $S
= (\zz^+)^{<\omega}$, and for each $m \geq 2$,
$$\ell^1(S,\,\omega_m) = {\cal U}_{m} \cong
\widehat{\bigvee}_{\omega_{m}}E,$$ where $\omega_m$ is a weight on
$\zz^+$ defined by $$\omega_m(|r|) = m^{|r|},\;r \in S =
(\zz^+)^{<\omega}$$ (or, one may take $\omega_m$ as a weight on
$S$ defined by $\omega_m(r) = m^{|r|}$). Then, $${\cal U} =
\bigcap_{m\in\nn}{\cal U}_m = \ell^1(S,\, W) \cong
\bigcap_{m\in\nn}\widehat{\bigvee}_{\omega_{m}}E =
\widehat{\bigvee}_{W}E.$$ All these algebras are graded
subalgebras of $\F_\infty$. The non-commutative analogues are the
Banach algebras $\ell^1_{nc}(S_{nc},\, \omega)\,=\,{\cal
U}_{1_{nc}}$, where $\omega\,\equiv\,1$ on $S_{nc}$, and for each
$m\,\geq\,2$, $$\ell^1_{nc}(S_{nc},\, \omega_m)\,=\,{\cal
U}_{m_{nc}}\,\cong\,\widehat{\bigotimes}_{\omega_{m}}E,$$ where
$\omega_m$ is a weight on $\zz^+$ defined by
$$\omega_m(|r|)\,=\,m^{|r|}, \;r\,\in\,S_{nc},
\;|r|\,=\,\textrm{rank}\,r$$ (or, one may take $\omega_m$ as a
weight on $S_{nc}$ defined by $\omega_m(r)\,=\,m^{|r|}$). Then,
$${\cal U}_{nc}\,=\,\bigcap_{m\in\nn}{\cal
U}_{nc_{m}}\,\cong\,\bigcap_{m\in\nn}\widehat{\bigotimes}_{\omega_{m}}E\,=\,\widehat{\bigotimes}_{W}E$$
~\cite[\S 2]{15}, ~\cite[\S 9, \S 10]{14} and ~\cite[\S 4]{30}.
All these algebras are graded subalgebras of $\cal B$. We remark
that the map $\epsilon$, restricted to ${\cal U}_m$, is an
isometric isomorphism of ${\cal U}_m$ onto
$\widehat{\bigvee}_{\omega_{m}}E$ (for $m\,=\,1$, see Remarks A.
1. above) and the same map $\epsilon$, restricted to $\cal U$, is
an isometric isomorphism of ${\cal U}$ onto
$\widehat{\bigvee}_{W}E$ ~\cite[p. 142]{14}.

It is important to note the following chains of (dense) continuous
inclusions of certain algebras and their corresponding character
spaces. We have
\begin{equation} \label{inclusion}
{\cal U}\,\hookrightarrow\,{\cal U}_{m}\,\hookrightarrow\,{\cal
U}_{1}\,=\,\ell^1(S,
\{1\})\,\cong\,\widehat{\bigvee}_{\{1\}}E\,\hookrightarrow\,\widehat{\bigvee}E,
\end{equation}
where the density follows from the fact that all algebras are
countably generated by the monomials $X_0, X_1,\dots$, and the
Arens-Michael representation of $\cal U$. By ~\cite[Lemma
3.2.5]{19}, we have
\begin{equation} \label{spectra}
M(\widehat{\bigvee}_{\{1\}}E)\,\cong\,\ell_{[1]}^\infty\,\cong\,M(\ell^1(S,\,
\{1\}))\,\hookrightarrow\,M({\cal
U}_m)\,\cong\,\ell_{[m]}^\infty\,\hookrightarrow\,M({\cal
U})\,\cong\,\ell^\infty
\end{equation} ($w^*$-topology), where the homeomorphisms are due to Theorem 2.1
and ~\cite{15}, respectively. We also have non-commutative
analogues of \eqref{inclusion} and \eqref{spectra}, but with the
same character spaces in \eqref{spectra} (e.g., $M({\cal
U}_{nc})\,\cong\,\ell^\infty\,\cong\,M(\cal U)$).

Next, we {\it claim} that the test case $\cal U$ (resp., ${\cal
U}_{nc}$) admits a Fr\'{e}chet topology $\tau$, if there exists a
discontinuous character on $\cal U$ (resp., on ${\cal U}_{nc}$).
The crucial point of the main theorem should be emphasized. It is
a surprising consequence of the fact that one is able to show that
the test case $\cal U$ (resp., ${\cal U}_{nc}$) is, indeed, a
weighted Fr\'{e}chet symmetric (tensor) algebra
$\widehat{\bigvee}_{W}E$ (resp., $\widehat{\bigotimes}_{W}E$) as
above. So if one starts with a discontinuous character on $\cal U$
(resp., on ${\cal U}_{nc}$), then one can generate a bizarre
Fr\'{e}chet topology on $\cal U$ (resp., on ${\cal U}_{nc}$) by
applying the Read's method in the following\\

\noindent{\bf Proof of the main theorem.}

We first show that if $\phi$ is a discontinuous character on $\cal
U$ (resp., on ${\cal U}_{nc}$), then
$\phi|_{\ell^1}\,=\,\lambda_1$ on $\ell^1$ is a discontinuous
linear functional on ${\cal U}^{(1)}\,=\,{\cal U}_{nc}^{(1)}$. We
recall that $\cal U$ and ${\cal U}_m$ ($m\,\in\,\nn$) are graded
subalgebras of $\F_\infty$; i.e.,
$${\cal U}\,=\,\sum_{p=0}^\infty{\cal U}^{(p)}\; \textrm{and}\; {\cal U}_{m}
= \sum_{p=0}^\infty{\cal U}_{m}^{(p)},\; \textrm{with}\; {{\cal
U}_{1}}^{(1)}\,=\,\ell^1(\zz^+)\,=\,\hat{\bigvee}^{1}E$$
~\cite{14}. So, ${\cal U}^{(p)}\,=\,\bigcap_{m\,\in\,\nn}{\cal
U}_{m}^{(p)}$ is the subspace of $p$-homogeneous formal power
series $\sum_{r\,\in\,(\zz^+)^{p}}\alpha_rX^r$, and, for each $m$
and $p$, there is a continuous, dense embedding from ${\cal
U}_{m}^{(p)}$ into ${\cal U}_{1}^{(p)}$, with ${\cal U}_{m}^{(p)}$
a closed linear subspace of ${\cal U}_{m}$ spanned by the
monomials $X^r$, $|r|\,=\,p$. So, clearly, for each $m$, ${\cal
U}_{m}^{(1)}$ is also isomorphic with ${\cal
U}_{1}^{(1)}\,=\,\ell^1(\zz^+)\; (\|\cdot\|_1\,\sim\,\|\cdot\|_m
\;\textrm{on} \;{\cal U}_{m}^{(1)},\; \textrm{where}\; {\cal
U}_{m}^{(1)}\,=\,\{\sum_{i=1}^\infty\alpha_iX_i\,:\,\sum_{i=1}^\infty|\alpha_i|m\,<\,\infty\}$).
In particular,
\begin{equation} \label{up} {\cal
U}^{(p)}=\bigcap_{m \in \nn}{\cal U}_{m}^{(p)}=\{\sum_{r \in
S}\alpha_rX^r \in {\cal U}_{1}:\sum_{r \in
S,\;|r|=p}\,|\alpha_r|m^p < \infty\;\textrm{for all}\,m \in \nn\}.
\end{equation}
Therefore,
\begin{equation} \label{U-only} {\cal U} =
\sum_{p=0}^\infty(\bigcap_{m\,\in\,\nn}{\cal U}_{m}^{(p)})\, =\,
\bigcap_{m\,\in\,\nn}(\sum_{p=0}^\infty{\cal U}_{m}^{(p)})\, =\,
\bigcap_{m\,\in\,\nn}{\cal U}_m,
\end{equation}
since ${\cal U}_m\,=\,\sum_{p=0}^\infty{\cal U}_{m}^{(p)}$. In
particular, \begin{equation} \label{U1} {\cal U}^{(1)}=\bigcap_{m
\in \nn}{\cal
U}_{m}^{(1)}=\{\sum_{i=1}^\infty\alpha_iX_i:\sum_{i=1}^\infty|\alpha_i|m
< \infty,\;\textrm{for all}\,m \in \nn\},\;{\cal U}_{m}^{(1)}
\cong \ell^1(\zz^+)
\end{equation}
as above. Here, \begin{equation} \label{U1-another} {\cal
U}_{1}=\sum_{p=0}^\infty{\cal
U}_{1}^{(p)}=\sum_{p=0}^\infty\ell^1(S^{(p)});\;{\cal
U}_{1}^{(1)}=\ell^1(S^{(1)}) \cong \ell^1(\zz^+) =
\hat{\bigvee}_{\{1\}}^{1}E.
\end{equation}
Also, ${\cal U}_{m+1}^{(p)}\,\subset\,{\cal U}_{m}^{(p)}$ for all
$m\,\in\,\nn$ and $p\,\in\,\zz^+$. We remark that ${\cal U}^{(1)}$
is a closed linear subspace of $\cal U$, spanned by $X_i$, whereas
${\cal U}_m^{(1)}$ is a closed linear subspace of ${\cal U}_{m}$
for each $m$, spanned by $X_i$ (the latter is a Banach space and
the former is a Fr\'{e}chet space). Obviously, all
$X_i\,\in\,{\cal U}^{(1)}$, but $\sum_{i=1}^\infty\frac{X_i}{i}$
does not belong to ${\cal U}^{(1)}$ (resp., does not belong to
even ${\cal U}_{1}^{(1)}$) whereas
$\sum_{i=1}^\infty\frac{X_i}{i^2}\,\in\,{\cal U}^{(1)}$ (resp.,
belongs to ${\cal U}_{1}^{(1)}$).

Assume that there is a discontinuous character $\phi$ on ${\cal
U}\,=\,\bigcap_{m\,\in\,\nn}{\cal U}_m$. Then, for each
$m\,\in\,\nn$, $\phi$ is a discontinuous character on the normed
algebra $({\cal U}, q_m)$, because if it is continuous for some
$m$, then there is $c\,>\,0$ such that
$$|\phi(f)|\,\leq\,cq_m(f),\; f\,\in\,({\cal U},\,(q_m)),$$ a
contradiction of the assumption that $\phi$ is discontinuous on
$\cal U$. Now, $\phi$ can be extended to a discontinuous linear
functional $\lambda_1$ on $({\cal U}_1, q_1)$ using a Hamel basis
of a complement of $\cal U$ in ${\cal U}_1$ (remark that there are
infinitely many ways of extending $\phi$; also, we cannot extend
$\phi$ as a discontinuous character on ${\cal U}_1$). (Similarly,
we can also extend $\phi$ to a discontinuous linear functional on
${\cal U}_m$ for each $m\,\in\,\nn$.) By Lemma 2.3,
$q_1\,\sim\,\|\cdot\|_{\ell^1}$ on $\ell^1$, and so, w.l.o.g.,
$\lambda_1$ (again denoted by $\phi$ below) is also a
discontinuous linear functional on
$\ell^1\,=\,\hat{\bigvee}_{\{1\}}^{1}E\,\cong\,{\cal
U}_{1}^{(1)}$. Alternatively, by a remark below ~\cite[Corollary
10]{8}, the existence of a discontinuous character on $\cal U$
induces a discontinuous linear functional on $\ell^1$. So,
$\phi|_{\ell^1} = \lambda_1$ is also a discontinuous linear
functional on ${\cal U}_{m}^{(1)}$ for each $m \in \nn$ (the
continuous embedding from ${\cal U}_{m}^{(1)}$ into ${\cal
U}_{1}^{(1)}$ actually turns out to be a topological isomorphism).
Thus, $\phi|_{\ell^1}$ is also a discontinuous linear functional
on ${\cal U}^{(1)}$ (because if it is a continuous linear
functional on ${\cal U}^{(1)} = \bigcap_{m \in \nn}{\cal
U}_{m}^{(1)}$, then, for each $m \in \nn$, it is also a continuous
linear functional on ${\cal U}_{m}^{(1)}$, a contradiction of the
fact that it is discontinuous on $\ell^1(\zz^+) \cong {\cal
U}_{1}^{(1)}$ due to our assumption).

Now, we choose $\phi(X_n) = 1$ for all $n \in \zz^+$, and using
the axiom of choice, we extend $\phi$ to a Hamel basis of ${\cal
U}^{(1)}$, so that $\phi$ becomes a discontinuous linear
functional on ${\cal U}^{(1)}$. This idea was actually used by
\.Zelazko in ~\cite{41} while defining a (discontinuous) linear
functional on $\textrm{Hol}(\cc)$, much before Read
~\cite[Definition 1.7]{33}. Read used it to define a discontinuous
linear functional on a closed linear subspace ${\cal A}^{(1)}$ of
$\cal B$, spanned by $X_i$. In our case, $${\cal U}^{(1)} \cong
\textrm{Hol}(\cc),\;\sum_n\alpha_nX_n \mapsto \sum_n\alpha_nz^n,$$
as shown in Remarks B. 2. below. So, $$H_0 = \{\psi \in
\textrm{Hol}(\cc) : \psi(z) = e^{z^n},\;n \in \nn\},$$ is a
countable linearly independent subset of $H$ ~\cite[Lemma 1]{41},
and $f$ could be defined on $\textrm{Hol}(\cc)$ using this $H_0$.
We shall generate another Fr\'{e}cet topology $\tau$ (inequivalent
to the usual one $\tau_0$, generated by the norms $q_m$) on ${\cal
U}_{nc}$ (and thus, on $\cal U$ as well), using this discontinuous
linear functional. Our argument here is kept short because it uses
the key ideas involved in the Read's method (so we follow
notations in align with Read's notations), but the Fr\'{e}chet
topologies on ${\cal U}_{nc}$ (resp., on $\cal U$) that we are
working with, are entirely different.

We start with a remark that the derivations on a
(non-)commutative, semi-simple Fr\'{e}chet algebra $({\cal U},
(q_m))$ (resp., $({\cal U}_{nc}, (q_m))$) are continuous ~\cite{7}
(Corollary 4.3 below). However, we shall show that the derivation
$\partial/\partial X_0$ is discontinuous w.r.t. $\tau$, a
contradiction. We shall also require in a future proof the
following ``locally finite" linear map
\begin{equation} \label{lf-maps} T\,:\,{\cal
U}_{nc}\,\rightarrow\,{\cal
U}_{nc},\;T(\sum_{r\,\in\,S_{nc}}\alpha_rX^{\otimes\,r})\,=\,\sum_{r\,\in\,S_{nc}}\alpha_rX_r
\end{equation}
(actually, this $T$ is the composition of the inclusion map with
the $T$, discussed by Read in ~\cite{33}; also, one requires to
restrict the range as well). Similarly we may define ``locally
finite" linear maps ${\cal U}_{nc} \rightarrow \cal U$, ${\cal U}
\rightarrow {\cal U}_{nc}$ and $\cal U \rightarrow \cal U$. They
are precisely the linear maps between these two spaces that are
continuous w.r.t. their usual Fr\'{e}chet topologies (two such
mappings were considered in the proof of Proposition 3.1 above).

We shall require the concept of ``tensor products by rows", taken
from ~\cite{33}. First, for each $n \in \zz^+$, let $P_n : {\cal
U}_{nc} \rightarrow {\cal U}_{nc}$ be the linear map such that
$P_n(1) = 0$, and\\
\begin{equation}
P_n(X_{i_{1}}\,\otimes\,X_{i_{2}}\,\otimes\,\cdots\,\otimes\,X_{i_{m}})
= \begin{cases} 0, ~&\text{when}\;i_{1}\,\neq\,n,\\
X_{i_{2}}\,\otimes\,\cdots\,\otimes\,X_{i_{m}}, ~&\text{when}\;
i_{1}\,=\,n. \end{cases}
\end{equation}
Thus $P_n$ takes the quotient on division from the left by $X_n$,
and discards the remainder. Let $\pi : {\cal U}_{nc} \rightarrow
{\cal U}$ be the natural map that the $X_i$ commute, i.e., the
locally finite map such that $\pi(X^{\otimes {\bf i}}) = X^{{\bf
i}}$ for all ${\bf i}$. Then $\pi|_{\tilde{\sigma}({\cal
U}_{nc})}$ is bijective and $\epsilon|_{{\cal U}} : {\cal U}
\rightarrow \tilde{\sigma}({\cal U}_{nc})$ be its inverse (recall
that $\epsilon : {\cal U} \rightarrow {\cal U}_{nc}$ is a
continuous monomorphism such that $\epsilon = \oplus \epsilon_p$,
where $\epsilon_p : {\cal U}^{(p)} \rightarrow {\cal
U}_{nc}^{(p)}$ is a continuous linear embedding for each $p$).
Thus $\epsilon$ is the natural right inverse to $\pi$.

Next, we have a linear functional $\phi : {\cal U}^{(1)}
\rightarrow \cc$ such that $\phi(X_n) = 1$ for all $n \in \zz^+$,
as discussed earlier. Recall that ${\cal U}_{nc}$ is a graded
subalgebra of $\cal B$. Write ${\cal U}_{nc}^{(p)}$ for the
subspace of $p$-homogeneous formal power series $\sum_{{\bf i}\,
\in\, S_{nc}^{(p)}}b_{{\bf i}}X^{\otimes \,{\bf i}}$ and write
${\cal U}^{(p)} = \pi({\cal U}_{nc}^{(p)}), \;{\cal U}_{nc}^{(1)}
= {\cal U}^{(1)}$. If ${\bf b} \in {\cal U}_{nc} =
\oplus_{p=0}^\infty{\cal U}_{nc}^{(p)}$ we write $({\bf
b}^{(p)})_{p\,\in\,\nn}$ such that ${\bf b} = \sum_{p=1}^\infty
{\bf b}^{(p)}$. If $\phi_1,\,\phi_2\,:\,{\cal
U}^{(1)}\,\rightarrow\,\cc$ are linear functionals, we define the
``tensor product by rows" $\phi_1 \otimes \phi_2 : {\cal
U}_{nc}^{(2)} \rightarrow \cc$ by
\begin{equation} \label{tensor-rows} \phi_1\,\otimes\,\phi_2\,({\bf
b})\,=\,\phi_1(\sum_{j=0}^\infty\,\phi_2(P_j{\bf b})X_j).
\end{equation}
Tensor product by rows of $n$ linear functionals are then defined
inductively by \begin{equation} \label{tensor-rows-several}
\otimes_{i=1}^n\,\phi_i({\bf
a})\,=\,\phi_1(\sum_{j=0}^\infty\,\otimes_{i=2}^n\,\phi_i(P_j{\bf
a})X_j).
\end{equation}
We see that ~\cite[Lemma 1.10]{33} holds true for the elements
${\bf a} \in {\cal U}_{nc}^{(r)},\;{\bf b} \in {\cal
U}_{nc}^{(n-r)}$ for $0 \leq r \leq n$, and functionals
$(\phi_i)_{i = 1}^n$ on ${\cal U}^{(1)}$. Then we have
\begin{lem} \label{Cor_Read} If $(\phi_n)_{n = 0}^\infty$ is a sequence of
linear functionals on ${\cal U}^{(1)}$, then, for each
$n\,\in\,\zz^+$, the seminorm $\|\cdot\|_n$ on ${\cal U}_{nc}$
given by
\begin{equation} \label{norms} \|{\bf a}\|_{n}\,=\,|{\bf
a}^{(0)}| \,+\,\sum_{r=1}^{\infty}\,\sum_{{\bf
i}\,\in\,(\zz^{+})^{r}}\,|\otimes_{j=1}^{r}\phi_{i_{j}}({\bf
a}^{(r)})|
\end{equation}
is a submultiplicative seminorm. $\hfill \Box$
\end{lem}
As discussed in ~\cite{33}, since the order of appearance of the
$\phi_{i_j}$ can be permuted arbitrarily in \eqref{norms}, one has
$\|\tilde{\sigma}\|_n\,=\,1$ for all $n$, where $\tilde{\sigma}$
is as in Proposition 3.1 above. Hence, the seminorms $\|\cdot\|_n$
are also submultiplicative seminorms on $\cal U$, when $\cal U$ is
identified with the linear subspace $\tilde{\sigma}({\cal U}_{nc})
\subset {\cal U}_{nc}$ (the multiplication of $\cal U$ is then
implemented by $({\bf a}, {\bf b}) \rightarrow \tilde{\sigma}({\bf
a} \otimes {\bf b})$). However, instead of the `usual' coordinate
linear functionals, discussed in ~\cite{33}, we, here, need the
`weighted' coordinate linear functionals, defined as follows, in
order to give Fr\'{e}chet topologies $\tau$ and $\tau_0$ on ${\cal
U}_{nc}$ below. Let
\begin{equation} \label{phi-nm} \phi_n\,:\,{\cal
U}^{(1)}\,\rightarrow\,\cc,\;\phi_n(\sum_{n=0}^\infty\alpha_nX_n)\,=\,\alpha_n\cdot\,n\;\;
\textrm{for each}\, n\,\in\,\nn,
\end{equation}
then $(\phi_n)_{n \in \nn}$ is a sequence of weighted coordinate
linear functionals on ${\cal U}^{(1)}$. Let $(\phi_n)_{n \in
\zz^+}$ be a sequence of linear functionals on ${\cal U}^{(1)}$ as
follows: (a) $\phi_0 = \phi$, the discontinuous linear functional
defined above; and (b) for $n \in \nn$, $\phi_n$ be the weighted
coordinate functional on ${\cal U}^{(1)}$. Apply Lemma
~\ref{Cor_Read} to $(\phi_n)$ to define a locally multiplicatively
convex topology $\tau$ on ${\cal U}_{nc}$.

We claim that $({\cal U}_{nc}, \tau)$ is a Fr\'{e}chet algebra.
Since $\tilde{\sigma}$ is a $\tau-$continuous projection, the
subspace $${\cal U} = \tilde{\sigma}({\cal U}_{nc}) =
\textrm{ker}\,(I - \tilde{\sigma})$$ is closed, so $({\cal U},
\tau)$ is a commutative Fr\'{e}chet algebra. Note that if we prove
that ${\cal U}_{nc}$ is, in fact, complete under the topology
$\tau$, then we have our desired result. For since such a
Fr\'{e}chet topology $\tau$ is inequivalent to the usual
F\'{e}chet topology $\tau_0$ on $\cal U$ (resp. on ${\cal
U}_{nc}$). We then arrive at a contradiction of the fact that
$({\cal U}, \tau_0)$ (resp., $({\cal U}_{nc}, \tau_0)$) is a
(non-)commutative, semi-simple Fr\'{e}chet algebra with a unique
Fr\'{e}chet topology by by ~\cite{6, 28, 29, 14} (resp., by
Corollary 4.3 below). Since $$\phi_0(X_N - X_0) =
0\;\textrm{for}\; N
> 0$$ one sees that $$\|X_N - X_0\|_n = 0\; \textrm{for all}\; n\; \textrm{with}\; n < N;$$
hence $X_N \rightarrow X_0$ in $\tau$, so one has $X_N - X_0
\rightarrow 0$. One may then use the derivation
$\partial/\partial\,X_0$ to see that $\partial/\partial\,X_0(X_0 -
X_n) = 1$, and so $1$ is in the separating subspace ${\cal
S}(\partial/\partial\,X_0)$, concluding that the separating
subspace is the whole algebra, be it $\cal U$ or ${\cal U}_{nc}$.
All that remains is to prove ${\cal U}_{nc}$ complete.

Let $\tau_0$ be the ``usual" topology that makes ${\cal U}_{nc}$ a
Fr\'{e}chet algebra (cf. ~\cite[\S2]{15} and ~\cite[\S 9]{14}).
One sees that $\tau_0$ could be obtained by applying Lemma 3.2 to
the sequence $(\phi_n)$, where $\phi_0$ is the usual continuous
coordinate functional, defined by
$\phi_0(\sum_{n=0}^\infty\alpha_nX_n) = \alpha_0$, and for $n \in
\nn$, $\phi_n$ is the weighted coordinate functional, defined as
above, inducing submultiplicative norms $|\cdot|_n$, where
\begin{equation} \label{norms-other} |{\bf a}|_{n}\,=\,|{\bf a}^{(0)}|
\,+\,\sum_{r=1}^{\infty}\,\sum_{{\bf
i}\,\in\,(\zz^{+})^{r}}\,|{\bf a}^{(r)}|n^{|{\bf i}|}.
\end{equation}
One may say that $\tau_0$ is the topology of ``weighted"
convergence w.r.t. all these norms $|\cdot|_n$; that is, $\tau_0$
is the topology of ``weighted" convergence of all the coefficients
$a_{{\bf i}}$.

Next, with an eye on (3.11), we define a linear map
$\Psi\,:\,{\cal U}_{nc}\,\rightarrow\,{\cal U}_{nc}$ by
\begin{equation} \label{linear map} \Psi({\bf
a})\,=\,a^{(0)}\,+\,\sum_{r=1}^{\infty}\,\sum_{{\bf
i}\,\in\,(\zz^{+})^{r}}\,X^{\otimes \,{\bf i}}\,\cdot
\,\otimes_{j=1}^{r}\phi_{i_{j}}({\bf a}^{(r)}).
\end{equation}
We note that $\Psi : ({\cal U}_{nc}, \tau) \rightarrow ({\cal
U}_{nc}, \tau_0)$ is continuous because convergence under $\tau$
is precisely convergence of all the weighted coordinate
functionals $\otimes_{j=1}^{r}\phi_{i_{j}}({\bf a}^{(r)})$ that
are involved in $\Psi({\bf a})$.
\begin{thm} \label{Thm 2.3_Read} $\Psi$ is bijective.
\end{thm}
\noindent{\bf Proof.} The proof is the same as that of
~\cite[Theorem 2.3]{33}, with a remark that one requires to
replace $\cal B$ by ${\cal U}_{nc}$ throughout that proof. The
readers notice that the Read's proof was purely algebraic, and so
one survives under the replacement, because ${\cal U}_{nc}$ is a
graded subalgebra of $\cal B$, so it inherits all the graded
algebraic structure that $\cal B$ has. We give a proof for the
sake of completeness.

$\Psi$ evidently maps ${\cal U}_{nc}^{(r)}$ to ${\cal
U}_{nc}^{(r)}$ for each $r\,\in\,\zz^+$. It is therefore enough to
show that $\Psi^{(r)}\,=\,\Psi|_{{\cal U}_{nc}^{(r)}}$ is a
bijection ${\cal U}_{nc}^{(r)}\,\rightarrow\,{\cal U}_{nc}^{(r)}$
for each $r$.

When $r\,=\,0$, we have the identity map on $\cc$, which is
bijective. For $r\,>\,0$, we proceed as follows.

\noindent {\bf Definition.} For $0 \leq k \leq r$, let $\nn^{(r,\,
k)}$ be the set of all ${\bf i} = (i_1, i_2, \dots, i_r) \in
(\zz^+)^r$ such that exactly $k$ of the $i_j$ are equal to zero.
Let ${\cal U}_{nc}^{(r,\, k)}\,\subset\,{\cal U}_{nc}$ be the
subspace consisting of all power series of the form $\sum_{{\bf i}
\in \nn^{(r,\, k)}}a_{\bf i}X^{\otimes \,{\bf i}}$.

Obviously, ${\cal U}_{nc}^{(r)}\,=\,\oplus_{k = 0}^r {\cal
U}_{nc}^{(r,\, k)}$; let $\mu_{(r, \,k)}$ be the projection onto
${\cal U}_{nc}^{(r,\, k)}$ parallel to the others. We claim that
for each $k$, $\mu_{(r,\,k)}\Psi^{(r)}|_{{\cal U}_{nc}^{(r,\,
k)}}$ is equal to the identity map on ${\cal U}_{nc}^{(r,\, k)}$;
and we claim that $\Psi^{(r)}$ maps ${\cal U}_{nc}^{(r,\, k)}$
into $\oplus_{l = k}^r {\cal U}_{nc}^{(r,\, l)}$. Thus we shall
show that the action of $\Psi^{(r)}$ is ``lower triangular" and
nonsingular, with respect to the decomposition ${\cal
U}_{nc}^{(r)}\,=\,\oplus_{k = 0}^r {\cal U}_{nc}^{(r,\, k)}$;
therefore, $\Psi^{(r)}$ - and hence $\Psi$ - is bijective.

In view of the definition of the map $\Psi$,
$\mu_{(r,\,k)}\Psi^{(r)}|_{{\cal U}_{nc}^{(r,\, k)}}$ is equal to
the identity map if and only if we have
\begin{equation} \label{2.1} \otimes_{i = 1}^{r}\phi_{l_{i}}({\bf
a})\,=\,a_{{\bf l}}|{\bf l}|
\end{equation}
for each ${\bf l}\,=\,(l_1,\dots,l_r)\,\in\,\nn^{(r,\,k)}$ and
${\bf a}\,\in\,{\cal U}_{nc}^{(r,\, k)}$; and $\Psi^{(r)}$ maps
${\cal U}_{nc}^{(r,\, k)}$ into $\oplus_{l = k}^r {\cal
U}_{nc}^{(r,\, l)}$ if and only if
\begin{equation} \label{2.2} \otimes_{i = 1}^{r}\phi_{l_{i}}({\bf
a})\,=\,0
\end{equation}
for all $m\,<\,k,\;{\bf l}\,=\,(l_1,\dots,
l_r)\,\in\,\nn^{(r,\,k)}$ and ${\bf a}\,\in\,{\cal U}_{nc}^{(r,\,
k)}$.

Let us prove (\ref{2.1}) and (\ref{2.2}) together, by induction on
$r$.

When $r\,=\,1$, statement (\ref{2.1}) demands that $\phi_l({\bf
a})\,=\,a_l$ if either $l\,=\,0$ and ${\bf a}\,\in\,{\cal
U}_{nc}^{(1,\, 1)}$, or $l\,>\,0$ and ${\bf a}\,\in\,{\cal
U}_{nc}^{(1,\, 0)}$. For $l\,>\,0$ the functional $\phi_l$ is the
weighted coordinate functional anyway, so the assertion is true
when $l\,>\,0$. But when $l\,=\,0$, we note that ${\cal
U}_{nc}^{(1,\, 1)}$ consists solely of multiples $a_0X_0$; since
$\phi_0(X_0)\,=\,1$, the assertion is true for $l\,=\,0$ also.
Statement (\ref{2.2}) demands that $\phi_l({\bf a})\,=\,0$ if
$l\,>\,0$ and ${\bf a}\,\in\,{\cal U}_{nc}^{(1,\, 1)}$. But in
that case, ${\bf a}$ is a multiple of $X_0$ and $\phi_l$ is the
weighted coordinate functional for $l\,>\,0$; so the assertion is
again correct.

When $r\,>\,1$, we proceed by induction. Equation (3.8) gives us
\begin{equation} \label{2.3} \otimes_{i = 1}^{r}\phi_{l_{i}}({\bf
a})\,=\,\phi_{l_1}(\sum_{j=0}^\infty\,\otimes_{i=2}^r\,\phi_{l_i}(P_j{\bf
a})X_j).
\end{equation}

To prove statement (\ref{2.1}), consider the case when ${\bf l}
\in \nn^{(r,\,k)}$ and ${\bf a} \in {\cal U}_{nc}^{(r,\, k)}$.
Note that $P_j{\bf a}$ lies in ${\cal U}_{nc}^{(r - 1,\, k)}$ if
$j > 0$, but in ${\cal U}_{nc}^{(r - 1,\, k - 1)}$ if $j = 0$,
because the division on the left by $X_0$ removes a factor $X_0$
from each monomial $X^{\otimes {\bf i}}$.

If $l_1\,=\,0$, then $(l_2,\dots,l_r)\,\in\,\nn^{(r - 1, \,k -
1)}$ so $\otimes_{i = 2}^r \phi_{l_i}$ will (by induction
hypothesis) annihilate any vector in ${\cal U}_{nc}^{(r - 1,\,
k)}$, but will send ${\bf b}\,\in\,{\cal U}_{nc}^{(r - 1,\, k -
1)}$ to $b_{l_2,\dots,\,l_r}|{\bf l}|$. Accordingly, when
$l_1\,=\,0$, (\ref{2.3}) is equal to $\phi_0((P_0{\bf
a})_{l_2,\dots, \,l_r}X_0)\,=\,(P_0{\bf a})_{l_2,\dots,
\,l_r}\,=\,a_{0,\, l_2,\dots,\,l_r}|{\bf l}|\,=\,a_{{\bf l}}|{\bf
l}|$ as required.

If on the other hand $l_1\,>\,0$, then $\phi_{l_1}$ is just the
weighted coordinate functional, so (\ref{2.3}) is equal to
$\otimes_{i=2}^r\,\phi_{l_i}(P_{l_1}{\bf a})$. In this case,
$(l_2,\dots,\,l_r)\,\in\,\nn^{(r - 1, \,k)}$ and $P_{l_1}{\bf
a}\,\in\,{\cal U}_{nc}^{(r - 1,\, k)}$, so by induction
hypothesis, (\ref{2.3}) is equal to $(P_{l_1}{\bf
a})_{l_2,\dots,\,l_r}\,=\,a_{{\bf l}}|{\bf l}|$ as before.

This deals with statement (\ref{2.1}). For (\ref{2.2}), we
consider ${\bf a}\,\in\,{\cal U}_{nc}^{(r,\, k)}$ as before, but
now use sequences ${\bf l}\,\in\,\nn^{r, \,m)}$ for some
$m\,<\,k$.

If $l_1\,=\,0$, then $(l_2,\dots,l_r)\,\in\,\nn^{(r - 1, \,m -
1)}$; since each $P_j{\bf a}$ is either in ${\cal U}_{nc}^{(r -
1,\, k)}$ or ${\cal U}_{nc}^{(r - 1,\, k - 1)}$, by induction
hypothesis the linear map $\otimes_{i = 2}^r \phi_{l_i}$
annihilates all of the vectors $P_j{\bf a}$, so (\ref{2.3}) is
equal to zero as required.

If $l_1\,>\,0$, then $\phi_{l_1}$ is the weighted coordinate
functional so (\ref{2.3}) is equal to
$\otimes_{i=2}^r\,\phi_{l_i}(P_{l_1}{\bf a})$. Since
$(l_2,\dots,l_r)\,\in\,\nn^{(r - 1, \,m)}$ and $P_{l_1}{\bf
a}\,\in\,{\cal U}_{nc}^{(r - 1,\, k)}$, by induction hypothesis
this is zero also. Thus (\ref{2.1}) and (\ref{2.2}) are proved,
and this completes the proof of Theorem 3.3. $\hfill \Box$

Let $(u_n)$ be a $\tau-$Cauchy sequence in ${\cal U}_{nc}$. Then
$(\Psi(u_n))$ is a $\tau_0-$Cauchy sequence in ${\cal U}_{nc}$.
Since $({\cal U}_{nc}, \,\tau_0)$ is complete, let $u$ be the
$\tau_0-$limit of the sequence $(\Psi(u_n))$. The $\tau-$limit of
the sequence $(u_n)$ is just $\Psi^{-1}(u)$. So $({\cal
U}_{nc},\,\tau)$ is complete under its seminorms $\|\cdot\|_n$.

The argument showing the rest of the assertions, given ${\cal
U}_{nc}$ complete, has already been given above. This completes
the proof of the main theorem. $\hfill \Box$

The analogue of our main theorem for the algebra $\cal B$ shows
that the non-commutative analogues of all avatars of the
Singer-Wermer conjecture do not hold in the Fr\'{e}chet case
~\cite[Remark 2.4]{31}.

\noindent {\bf Proof of Corollary.} By the main theorem, $\cal U$
(resp., ${\cal U}_{nc}$) admits two inequivalent Fr\'{e}chet
topologies, a contradiction, as discussed above. Also, by the main
theorem, the derivation $\partial/\partial\,X_0$ is discontinuous
on $\cal U$ (resp., on ${\cal U}_{nc}$), a contradiction as we
remark earlier. Since all characters of $\cal U$ are continuous,
by Proposition 3.1, all characters of ${\cal U}_{nc}$ are also
continuous. By ~\cite[Proposition 2.5]{18} and ~\cite[Remark
2.2]{15}, they are, indeed, the point evaluation mappings, that
is, if $\phi$ is a character on ${\cal U}_{nc}$, then there exists
$z\,\in\,\ell^\infty$ such that $\phi\,=\,\phi_z$, where
$\phi_z(f)\,=\,f(z)$ for all $f\,\in\,{\cal U}_{nc}$. $\hfill
\Box$
\begin{cor} \label{Patel_Michael} Every Fr\'{e}chet algebra (commutative
or not) is functionally continuous.
\end{cor}
\noindent{\bf Proof.} This is clear from Proposition 2.1,
~\cite[Remark 2.2]{15} and Corollary from \S 1. $\hfill \Box$

We have an interesting consequence as follows. First, recall that
Michael also asked whether every character on a commutative,
complete LMC-algebra is bounded ~\cite[\S 12, Question 2]{27}.
Dixon and Fremlin showed that the two problems are, in fact,
equivalent ~\cite{16}. Akkar gave in a nice interpretation of this
fact: the two problems are equivalent because every complete LMC
algebra is, as a ``bornological algebra", isomorphic to inductive
limit of a family of Fr\'{e}chet algebras ~\cite{1}. The
equivalence of the two questions also follows from the fact that
there is a complete LMC-algebra, produced by Craw ~\cite{9}, such
that the existence of a discontinuous character on some
commutative, unital Fr\'{e}chet algebra would imply the existence
of an unbounded character on the complete LMC-algebra. We have the
following
\begin{cor} \label{Michael_LMC} Every
complete LMC-algebra (commutative or not) is functionally bounded.
\end{cor}
\noindent {\bf Proof.} The technique of ~\cite{16} can easily be
applied to the non-commutative case as well (see the proof of
~\cite[Theorem 1]{16}). $\hfill \Box$\\

\noindent {\bf Remarks B.} 1. Until now, we know that $M(A)$ is
dense in $S(A)$, where $S(A)$ (resp., $M(A)$) is the space of all
(continuous) characters on $A$ w.r.t. the Gel'fand topology
~\cite[Lemma 10.1.1]{19}. By Corollary 3.4, we have $M(A) = S(A)$.

\noindent 2. Vogt studied the tensor algebras over a Fr\'{e}chet
space in ~\cite{39, 40}, mostly linear topological properties
(and, in particular, linear isomorphisms) of these algebras. A
particular topic of interest, close to the algebra $\cal U$ of all
entire functions on $\ell^\infty$, should be the study of the
space of entire functions on a nuclear sequence space, which is
the symmetric tensor algebra of a suitable nuclear reflexive
K$\ddot{o}$the space.

\noindent {\bf Fr\'{e}chet analogues of Sections 2 and 3, an
alternative of our approach.} It is clear that ${\cal U}^{(1)}$ is
isometrically isomorphic to the Beurling-Fr\'{e}chet algebra
$\ell^1(\zz^+, W)$, where $W = (\omega_m)$ is an increasing
sequence of weights on $\zz^+$ defined by $\omega_m(n) = m$ for
all $n, m \in \nn$ and $\omega_m(0) = 1$ for all $m \in \nn$
~\cite[Example 1.2]{4}. This is, in turn, isomorphic to the
Fr\'{e}chet algebra $\textrm{Hol}(\cc)$ of entire functions
~\cite[Example 1.4]{4}, which is the nuclear power series space of
infinite type ~\cite[\S 2]{40}. Note that the linear isomorphism
$\sum_{n=0}^\infty\alpha_nX_n \mapsto
\sum_{n=0}^\infty\alpha_nX^n$ between ${\cal U}^{(1)}$ and
$\ell^1(\zz^+, W)$ can be used to define a multiplication on the
Fr\'{e}chet space ${\cal U}^{(1)}$ which makes it a Fr\'{e}chet
algebra. Since $E = {\cal U}^{(1)}$ is a nuclear Fr\'{e}chet
space, it has an approximation property. If we apply fairly
straightforward arguments of Vogt ~\cite[\S 1]{40}, then it is
easy to see that the tensor algebra $T(E)$ (resp., the symmetric
algebra $S(E)$ in the commutative case) is, indeed, the algebra
${\cal U}_{nc}$ (resp., ${\cal U}$). In fact, one develops the
theory of topological tensor algebra over the Fr\'{e}chet space $E
= {\cal U}^{(1)}$ as discussed in \S 2. For example, write
$\hat{\bigotimes}^pE$ for the completion of $\bigotimes E$ w.r.t.
the projective tensor product metric induced by the metric $d$ on
$E$, then we have $\hat{\bigotimes}^pE = {\cal U}_{nc}^{(p)}$, the
subspace of $p$-homogeneous formal power series, and
$\widehat{\bigotimes}E = {\cal U}_{nc}$ the weighted Fr\'{e}chet
tensor algebra. Similarly, we have $\hat{\bigvee}^pE = {\cal
U}^{(p)} = SE(p)$, the closed subspace of $\hat{\bigotimes}^pE =
{\cal U}_{nc}^{(p)}$ and $\widehat{\bigvee}E = {\cal U}$ the
weighted Fr\'{e}chet symmetric algebra.

The tensor algebra $T(E) = {\cal U}_{nc}$ is a nuclear power
series space of infinite type by ~\cite[Lemma 2.1 and Theorem
3.1]{40}, which is further linearly isomorphic to the space $s$ of
all rapidly decreasing sequences by ~\cite[Theorem 4.1]{40}. Then,
this isomorphism equips the space $s$ with a multiplication which
turns it into a Fr\'{e}chet algebra, called $s_\bullet$. By
Corollary 3.4 above, the algebra $T(s) = s_\bullet = T(E) = {\cal
U}_{nc}$ is now functionally continuous. Thus, the algebra
$s_\bullet$ was, indeed, one of the test cases for Michael's
problem (cf. ~\cite[p. 190]{40}). However, our approach to the
algebras ${\cal U}_{nc}$ and ${\cal U}$ seem to be more convenient
for our purpose in this section.

\noindent 3. Read used his method of constructing another
Fr\'{e}chet topology on $\F_\infty$, inequivalent to the usual
Fr\'{e}chet topology on $\F_\infty$, to show that the famous
Singer-Wermer conjecture fails in the Fr\'{e}chet case whereas we
use the same method to affirmatively answer Michael's problem. As
a welcome bonus, we also establish some significant, long-standing
problems in automatic continuity theory in \S 4 below.

In an another paper, we exclusively discuss the first approach,
based on the Dales-McClure's and Dales-author-Read's methods of
generating a natural, totally discontinuous higher point
derivations $d_n$ at a discontinuous character $d_0 = \phi$ on the
test case ${\cal U}$, to affirmatively solve Michael's problems in
the commutative case. We remark that Dales and McClure used the
topological version of the tensor algebra over a Banach space to
obtain a somewhat negative result in automatic continuity theory
whereas we use the same technique to answer Michael's problem.

Thus these two approaches confirm the impression that one may not
expect the phenomena to hold in the Fr\'{e}chet algebra theory
that hold in the Banach algebra theory. However, we fortunately
answer Michael's problem in our favor ~\cite[Conclusion]{33}.
Thus, the situation on Fr\'{e}chet algebras is markedly different
from that on Banach algebras, and that a structure theory for
Fr\'{e}chet algebras behaves in a very distinctive manner from
Banach algebras.

This clearly establishes the fact that how all this works have
turned out to be the stepping stones for the complete solution to
Michael's problem since 2004, as discussed in Introduction.

\noindent 4. Allan exhibited a discontinuous homomorphism between
two commutative, unital Fr\'{e}chet algebras having certain
properties ~\cite[Theorem 8]{2}. In the construction of
discontinuous homomorphism, he used a continuous homomorphism from
$A$ into $\F$, induced by a natural, continuous higher point
derivation $(d_n)$ of infinite order on the Fr\'{e}chet algebra
$\textrm{Hol}(\cc)$ of entire functions at a continuous character
$d_0$. He used only the ZFC axioms (and not the CH). Similarly,
Read used the ZFC axioms alone to give an inequivalent Fr\'{e}chet
algebra topology on the algebra ${\cal F}_\infty$, in order to
show that the Singer-Wermer conjecture fails in the Fr\'{e}chet
case ~\cite{33}. Thus we, here, establish $NDH_{F,\,\cc}$, where
$F$ is any Fr\'{e}chet algebra (see ~\cite[\S 9]{10} for details).
In ~\cite{14}, Dales, author and Read also used only the ZFC
axioms to establish several important results in automatic
continuity theory, including an affirmative answer to the
Dales-McClure problem, and showing that the test algebra $\cal U$
for Michael's problem is, indeed, a FrAPS in $\F$. We remark that
Michael's problem arises in the ZFC.

\noindent 5. We have throughout assumed that our algebras have an
identity. It is a triviality to show that if a commutative
Fr\'{e}chet algebra has a discontinuous character, then so does
the algebra obtained by adjoining an identity
(see~\cite[Proposition 2.1]{15} or ~\cite[Remarks 4.5 (2)]{9}).
Thus, as far as Michael's problem is concerned, no loss of
generality is involved in our assumption.
\section{Applications to automatic continuity theory}
\noindent {\bf Some remarks in automatic continuity theory.} As
promised, we have affirmatively answered Michael's problem within
automatic continuity theory. In the past, there were some
significant works in this theory, which relate to our past (see
Introduction) and present works in some or the other way and
giving applications within this theory or in the theory of
(commutative) rings and algebras, and in the theory of SCV (e.g.,
Gleason's problem (1964) ~\cite{32}) as follows. Johnson
established the uniqueness of the complete norm for semi-simple
Banach algebras ~\cite{21}, and for BAPS ~\cite{22}. This result
was extended by Loy to FrAPS satisfying the equicontinuity
condition (E) ~\cite{25}. Thomas showed that the image of a
derivation on a Banach algebra is contained in the radical,
solving the Singer-Wermer conjecture affirmatively ~\cite{36}. The
non-commutative analogue of this conjecture is still open in the
Banach case ~\cite{37}. However, Read showed that this conjecture
fails in the Fr\'{e}chet case ~\cite{33}. Loy gave a simple method
to construct commutative Banach algebras with inequivalent
complete norms by using the discontinuity of derivations
~\cite{26}; see, also, ~\cite{24}.\\

\noindent {\bf Remarks C.} 1. Our first remark is on Michael's two
problems for commutative Fr\'{e}chet algebras ~\cite[\S 12]{27}.
Their non-commutative analogues also exist as discussed above, and
Dixon and Esterle discussed the test case ${\cal U}_{nc}$ in
~\cite{15}. Dales posed a question about the continuity of
characters on a commutative locally convex $(F)$-algebras
~\cite[Question 3 (ii)]{10}. Dixon discussed the non-equivalence
of Michael's two problems for more general commutative locally
convex $(F)$-algebras and commutative, complete locally convex
algebras (the reason is that he could not extend ~\cite[Theorem
1]{16} in this case). We shall show that every character on a
(non-)commutative $(F)$-algebra is automatically continuous
(Corollary 4.2 below). Thus, every non-commutative, complete
LMC-algebra is functionally bounded by Corollary 3.5, but this is
not true for every commutative, complete locally convex algebra.

\noindent 2. It is obvious that a discontinuous character on
$({\cal U},\,\tau_0)$ (resp., on $({\cal U}_{nc},\,\tau_0)$) is
also discontinuous on $({\cal U},\,\tau)$ (resp., $({\cal
U}_{nc},\,\tau)$) by ~\cite[Theorem 2.7]{18} (resp., by
Proposition 3.1 above). So, Corollary 6 of ~\cite{6}, stated for
commutative, semi-simple Fr\'{e}chet algebras, is, here,
established for {\it all} semi-simple Fr\'{e}chet algebras in a
more general form by Corollary 3.4 and Corollary 4.3 below.

\noindent 3. Another topology $\tau$ on $\cal U$ (resp., on ${\cal
U}_{nc}$) also contradicts the fact that every derivation is
automatically continuous on a (non-)commutative, semi-simple
Fr\'{e}chet algebra ~\cite{7} (Corollary 4.3 below), because, by
the main theorem, the derivation $\partial/\partial\,X_0$ is
discontinuous w.r.t. $\tau$.

\noindent 4. Carpenter ~\cite{6} showed that every commutative
semi-simple Fr\'{e}chet algebra $A$ admits a unique Fr\'{e}chet
topology. In the second paper, we have obtained another equivalent
form of Michael's problem ~\cite[p. 143]{10} by extending the
\v{S}ilov's theorem ~\cite[Theorem 2.3.3]{11}. What about the
non-commutative analogues of these results? That is, whether a
non-commutative, semi-simple Fr\'{e}chet algebra has a unique
Fr\'{e}chet topology ~\cite[Question 9]{10}. More generally,
whether every homomorphism from a Fr\'{e}chet algebra onto a
semi-simple Fr\'{e}chet algebra is automatically continuous (a
weaker version of an equivalent form of Michael's problem). There
is a parallel question about the continuity of derivations on a
non-commutative, semi-simple Fr\'{e}chet algebra. As far as we
know, the first problem was discussed by Johnson in the Banach
case ~\cite{21}. A generalization of the Johnson's result was
given by Jewell and Sinclair ~\cite{20}, which turned out to be a
starting point for Esterle and Thomas to generalize the Jewell and
Sinclair stability lemma to a more general $(F)$-space case
~\cite{17, 35} as well as for author to discuss the three
questions here. Also, we do not know any progress on these
questions since 1967, 1978 and 1976, resp.

It is interesting to note that we can affirmatively answer
Michael's problems as an application of an equivalent form of
Michael's problem above by taking $A = \cc$. In fact, if $\theta$
is a homomorphism from a Fr\'{e}chet algebra $B$ (commutative or
not) into a commutative, semi-simple Fr\'{e}chet algebra $A$, then
it is automatically continuous if and only if every character on a
Fr\'{e}chet algebra $B$ is continuous. However, we, here, take an
alternate approach which also answers the third question. For
this, we first remark that a generalization of ~\cite[Theorem
2]{20} holds true in the Fr\'{e}chet case as follows (also, we can
take a (non-)commutative $(F)-$algebra in the domain of an
epimorphism, in order to state the below theorem in a fuller
generality).
\begin{thm} \label{JS_Patel} Let $(B, \,(q_k))$ be a Fr\'{e}chet algebra such
that
\begin{enumerate}
\item[{\rm (i)}] for each infinite-dimensional closed (two-sided)
ideal $J$ in $B$ there is a sequence $(b_n)$ in $B$ such that the
sequence $(\overline{(Jb_n \dots b_1)})_{n \in \nn}$ of closed
ideals is strictly decreasing in $B$ (i.e., this sequence does not
stabilize); \item[{\rm (ii)}] $B$ contains no non-zero
finite-dimensional nilpotent two-sided ideals.
\end{enumerate} Then a homomorphism from a Fr\'{e}chet algebra onto $B$ is
continuous. In particular, $B$ has a unique Fr\'{e}chet algebra
topology. Also, every derivation on $B$ is continuous.
\end{thm}
\noindent{\bf Proof.} Let $\theta$ be a homomorphism from a
Fr\'{e}chet algebra $A$ onto $B$. Then, by ~\cite[\S 6, pp.
132-133]{14}, the separating space ${\cal S}(\theta)$ is a closed
(two-sided) ideal in $B$, and if it is infinite-dimensional, then
there is a sequence $(b_n) \in B$ such that $(\overline{({\cal
S}(\theta)b_n\dots b_1)})$ is a strictly decreasing sequence of
closed (two-sided) ideals in $B$ (i.e., this sequence does not
stabilize in $B$). Since $\theta$ is an epimorphism, there is a
sequence $(a_n) \in A$ such that $\theta(a_n) = b_n$ for all $n
\in \nn$. We obtain a contradiction by applying a common extension
of ~\cite[Lemma 1]{20} and ~\cite[Lemma 1.1a]{35} with $$R_n(y) =
yb_n\;(y \in B)\; \textrm{and}\; T_n(x) = xa_n\;(x \in A)$$ as
follows.

First, for a fixed sequence $(b_n)_{n\,\in\,\nn}$ in $B$,
$(\overline{({\cal S}(\theta)b_n\dots b_1)})$ is a sequence of
closed (two-sided) ideals in $B$. So, by the Arens-Michael
representation of $\overline{({\cal S}(\theta)b_n\dots b_1)}$, for
each $n\,\in\,\nn$, $$\overline{({\cal S}(\theta)b_n\dots
b_1)}\,=\,\lim\limits_{\longleftarrow} \overline{({\cal
S}(\theta)b_n\dots b_1)}^k;$$ note that, for each $k,\,
n\,\in\,\nn$, $$\overline{Q_k(\overline{({\cal S}(\theta)b_n\dots
b_1)})}^k\,=\,\overline{({\cal S}(\theta)b_n\dots b_1)}^k,$$ where
$Q_k\,:\,B\,\rightarrow\,B_k$ is a continuous algebra homomorphism
(cf. ~\cite[Lemma 2.1]{28}). Hence, for each $k,\, n\,\in\,\nn$,
$\overline{({\cal S}(\theta)b_n\dots b_1)}^k$ is a closed
(two-sided) ideal in $B_k$. Then, following ~\cite[Lemma
1.1a]{35}, given $k\,\in\,\nn$, there exists $n(k)\,\in\,\nn$ such
that
$$\overline{({\cal S}(\theta)b_m\dots b_1)}^k\,=\,\overline{({\cal
S}(\theta)b_{n(k)}\dots b_1)}^k\; \textrm{for all}\;
m\,\geq\,n(k)$$ (i.e., by ~\cite[\S 6, p. 133]{14} and the
Arens-Michael isomorphism theorem, for each $k\,\in\,\nn$, ${\cal
S}(\theta)^k$ is a {\it separating ideal} in $B_k$). Hence, by
passing to a suitable subsequence $(q_{n(k)})$ of $(q_k)$ defining
the same Fr\'{e}chet topology of $B$, we conclude, w.l.o.g., that
for each $k\,\in\,\nn$, there exits $n\,\in\,\nn$ such that
$$\overline{({\cal S}(\theta)b_m\dots b_1)}^k\,=\,\overline{({\cal
S}(\theta)b_n\dots b_1)}^k\; \textrm{for all}\; m\,\geq\,n.$$ So,
for such an $n$,
$$\lim\limits_{\longleftarrow}\overline{({\cal S}(\theta)b_m\dots
b_1)}^k\,=\,\lim\limits_{\longleftarrow}\overline{({\cal
S}(\theta)b_n\dots b_1)}^k\; \textrm{for all}\; m\,\geq\,n.$$
Hence, there exits $n\,\in\,\nn$ such that $$\overline{({\cal
S}(\theta)b_m\dots b_1)}\,=\,\overline{({\cal S}(\theta)b_n\dots
b_1)}\; \textrm{for all}\; m\,\geq\,n.$$ Conversely, suppose that
there exists $n\,\in\,\nn$ such that $$\overline{({\cal
S}(\theta)b_m\dots b_1)}\,=\,\overline{({\cal S}(\theta)b_n\dots
b_1)}\; \textrm{for all}\; m\,\geq\,n.$$ Then, by the
Arens-Michael isomorphism theorem, for such an $n$,
$$\lim\limits_{\longleftarrow}\overline{({\cal S}(\theta)b_m\dots
b_1)}^k\,=\,\lim\limits_{\longleftarrow}\overline{({\cal
S}(\theta)b_n\dots b_1)}^k\; \textrm{for all}\; m\,\geq\,n.$$ So,
by analogous arguments from ~\cite[(3.3.3), pp. 86-87]{19}, for
each $k\,\in\,\nn$, there exists $n(k)\,\in\,\nn$ such that
$$\overline{({\cal S}(\theta)b_m\dots b_1)}^k\,=\,\overline{({\cal
S}(\theta)b_{n(k)}\dots b_1)}^k\; \textrm{for all}\;
m\,\geq\,n(k).$$

Hence ${\cal S}(\theta)$ is a closed finite-dimensional ideal in
$B$. Each element $b$ of ${\cal S}(\theta)$ has connected spectrum
containing $0$ since $\sigma(b)\,=\,\bigcup_{k = 1}^\infty
\sigma(b_k)$ and, by ~\cite[Proposition 5.1.3 (ii)]{11}, each
element $b_k$ of ${\cal S}(Q_k\,\circ\,\theta)$ has connected
spectrum containing $0$. So each $b$ of ${\cal S}(\theta)$ is
nilpotent because ${\cal S}(\theta)$ is finite-dimensional. Thus
${\cal S}(\theta)$ is a nilpotent ideal, so ${\cal
S}(\theta)\,=\,\{0\}$, and $\theta$ is continuous. Alternatively
the proof from ${\cal S}(\theta)$ finite-dimensional to ${\cal
S}(\theta)\,=\,\{0\}$ may be handled as for derivations, below.

Let $D$ be a derivation from $B$ into $B$. The separating space
${\cal S}(D)$ is a closed two-sided ideal in $B$ ~\cite[\S 0]{33}.
Similar working to that for homomorphisms with
$$R_n(y)\,=\,T_n(y)\,=\,yb_n\;\textrm{for all}\; n\,\in\,\nn$$ in the
application of the above common extension of ~\cite[Lemma 1]{20}
and ~\cite[Lemma 1.1a]{35} and the formula $$(DT_n -
R_nD)(y)\,=\,D(yb_n) - D(y)b_n\,=\,yD(b_n)$$ show that ${\cal
S}(D)$ is finite-dimensional. Suppose that $${\cal
S}(D)\,\neq\,\{0\}.$$ The (Jacobson) radical of an ideal is the
intersection of the ideal and the radical of the algebra, and so
is an ideal in the algebra. Hence the radical of ${\cal S}(D)$ is
a finite-dimensional nilpotent ideal in $B$, and so is zero. The
Wedderburn structure theorem for finite-dimensional semi-simple
algebras implies that ${\cal S}(D)$ has an identity $e$. The
identity in an ideal in an algebra is a central idempotent in the
algebra. Thus
$$D(e)\,=\,D(e^2)\,=\,eD(e)\,+\,D(e)e\,=\,2eD(e)\,\Rightarrow\,
D(e)\,=\,2eD(e)\,=\,0.$$ This shows that $$D(eB)\,\subseteq\,eB\;
\textrm{and}\; D((1 - e)B)\,\subseteq\,(1 - e)B.$$ Since $eB$ is
finite-dimensional, the restriction of $D$ to $eB$ is continuous.
From the decomposition $B \,=\,(1 -e)B\,\oplus\,eB$ it follows
that ${\cal S}(D)$ is contained in $(1 - e)B$. Since $e$ is in
${\cal S}(D)$ and is not in $(1 - e)B$ this is a contradiction.
This completes the proof of theorem. $\hfill \Box$

A closed ideal $I$ in a Fr\'{e}chet algebra $A =
\lim\limits_{\longleftarrow}A_k$ is locally nilpotent if and only
if $I_k$ is a nilpotent ideal in $A_k$ for each $k$ (cf. ~\cite[p.
277]{2}). We now emphasize that no improvement in the condition
(ii) such as ``$B$ contains no non-zero finite-dimensional locally
nilpotent two-sided ideals", is possible here. The reason is that
since ${\cal S}(\theta)$ is a closed, finite-dimensional ideal,
every locally nilpotent element of ${\cal S}(\theta)$ is, indeed,
nilpotent, and so, for finite-dimensional two-sided ideals, the
notions of locally nilpotent ideals and nilpotent ideals coincide.

We remark that if $B$ is commutative in the above theorem, then
the condition (i) above can be expressed in a more neater form as
follows. For each infinite-dimensional closed ideal $J$ in $B$
there exists $b \in B$ such that the infinite-dimensional closed
ideal $\overline{(Jb)}$ is a proper subset of $J$.

\noindent {\bf Third approach to Michael's problem.}

It is a surprising consequence of Theorem 4.1 (with $B$ a Banach
algebra) that we have another approach to affirmatively answer
Michael's problems in the following corollary (see a comment on
~\cite[p. 144]{10}, or ~\cite{17}). Esterle remarked that
Corollary 3 (or Theorem 1) of ~\cite{17} cannot be applied in the
case $B = \cc$ for obvious reason. However, $B = \cc$ surely
satisfies both the conditions of ~\cite[Theorem 2]{20} (the first
condition is vacuously satisfied, and since $\cc$ is a field, the
second condition is also satisfied). In particular, $B$ does
satisfy ~\cite[Corollary 3]{17}, if the second condition is
replaced by ``if $B$ has no non-zero nilpotent ideal of finite
codimension" (equivalently, $B$ has no non-zero finite-dimensional
nilpotent ideal, which is the condition (ii) of ~\cite[Theorem
2]{20}). We can take $A$ a (non-)commutative $(F)-$algebra in the
below corollary, as mentioned in Remarks C. 1 above; also,
Esterle's stability lemma ~\cite[Lemma 1]{17} would suffice to
obtain a contradiction in the proof since $B$ is a Banach algebra.
\begin{cor} \label{Esterle_Patel} Let $A$ be a (non-)commutative
Fr\'{e}chet algebra, and $B$ be a Banach algebra as in Theorem 2
of ~\cite{20}. Then every epimorphism from $A$ onto $B$ is
continuous. In particular, every character on $A$ is continuous.
$\hfill \Box$
\end{cor}
\noindent {\bf Non-commutative analogues of some important results
in automatic continuity theory.}

Next, if $B$ is a non-commutative, semi-simple Fr\'{e}chet
algebra, then it has no nilpotent two-sided ideals. To see whether
the condition (i) of Theorem 4.1 holds, we follow the proof of
~\cite[Corollary 9]{20} in the Fr\'{e}chet case by working with an
infinite-dimensional irreducible left $B-$module $X$, which is a
Fr\'{e}chet module under the Fr\'{e}chet space topology. Then we
have the following
\begin{cor} \label{JS_Patel_Cor} Let $B$ be a non-commutative, semi-simple
Fr\'{e}chet algebra. Then a homomorphism from a Fr\'{e}chet
algebra onto $B$ is continuous. In particular, $B$ has a unique
Fr\'{e}chet algebra topology. Also, every derivation on $B$ is
continuous.
\end{cor}
\noindent {\bf Proof.} A semi-simple Fr\'{e}chet algebra contains
no non-zero nilpotent ideals. Let $J$ be a closed
infinite-dimensional ideal in a semi-simple Fr\'{e}chet algebra
$B$.

Suppose that there is an infinite-dimensional irreducible left
$B$-module $X$ such that $JX \neq \{0\}$. Let $(x_n)$ be a
linearly independent sequence in $X$. By induction using the
Fr\'{e}chet analogue of ~\cite[Lemma 7]{20} we choose a sequence
$(b_n)$ from $B$ so that $b_nb_{n - 1} \cdots b_1x_n = 0$ and
$\{b_n(b_{n - 1} \cdots b_1x_j) : j = n + 1,\, n + 2,\dots\}$ is
linearly independent. Then $\overline{(Jb_n \cdots b_1)}x_{n + 1}
\supseteq J(b_n \cdots b_1x_{n + 1})$ is non-zero because $Jy$ is
non-zero for each non-zero $y$ in $X$. An irreducible left module
over a Fr\'{e}chet algebra may be given a Fr\'{e}chet space
topology in which it becomes a Fr\'{e}chet module. With this
Fr\'{e}chet space topology on $X$ we obtain $$\overline{(Jb_{n +
1} \cdots b_1)}x_{n + 1}\,\subseteq\,\overline{J(b_{n + 1} \cdots
b_1x_{n + 1})}\,=\,\{0\}$$ so that $\overline{(Jb_n \cdots
b_1)}\,\supset\,\overline{J(b_{n + 1} \cdots b_1)}$ for all
$n\,\in\,\nn$.

We now suppose that $JX = \{0\}$ for each infinite-dimensional
irreducible left $B$-module $X$. Then there is an infinite
sequence $(X_n)$ of pairwise non-equivalent finite-dimensional
algebraically irreducible left $B$-modules such that $JX_n \neq
\{0\}$ for each $n \in \nn$. If there were no such sequence, then
$J$ would be finite-dimensional as it would be contained in the
sum of a finite number of finite-dimensional matrix algebras so
there would be a contradiction. By the Fr\'{e}chet analogue of
~\cite[Lemma 8]{20}, we choose a sequence $(b_n)$ from $B$ so that
$$b_nX_n = 0\; \textrm{and}\; b_nX_j = X_j\; \textrm{for}\; j\,=\,n + 1, \,n + 2,\dots.$$
Then $$\overline{(Jb_n \dots b_1)}X_{n + 1}\,=\,X_{n + 1}\;
\textrm{and}\; \overline{(Jb_{n + 1} \cdots b_1)}X_{n +
1}\,=\,\{0\}$$ so that $$\overline{(Jb_n \cdots
b_1)}\,\supset\,\overline{J(b_{n + 1} \cdots b_1)}\; \textrm{for
all}\; n\,\in\,\nn.$$ This completes the proof of Corollary 4.3.
$\hfill \Box$

We remark that Corollary 4.3 affirmatively answers ~\cite[Question
9]{10}. Esterle deduced in ~\cite{17} that, if $B$ is a Banach
algebra satisfying the condition (i) of Theorem 4.1, and if $B$
contains no non-zero finite-dimensional two-sided ideal (note the
variation from condition (ii) of Theorem 4.1), then every
epimorphism from an $(F)-$algebra onto $B$ is continuous. We see
that a generalization of Esterle's result does not hold in the
Fr\'{e}chet case. For example, the algebra $\F$ clearly satisfies
the conditions of Esterle's theorem, but, by ~\cite[Theorem
11.2]{14}, every homomorphism from an $(F)-$algebra onto $\F$ is
discontinuous, answering a Dales-McClure's problem in a stronger
form for a more general case (see ~\cite[Theorem 2.3]{13} for the
case the domain algebra a Banach algebra).\\

\noindent {\bf Acknowledgements.} This paper is based on a poster
presented at the Short Communications Satellite (SCS 2022, a part
of Virtual ICM 2022). The author thanks Professors Terence Tao and
Gerhard Paseman, Organizers of the SCS 2022, for providing me this
opportunity (my poster received the first rank among all entries).
We also take this opportunity to remember and thank significant
mathematicians, who contributed to Michael's problem in the past
and who, unfortunately, passed away before this work comes in
existence, especially, Stanis\l aw Mazur (1981), Richard Arens
(2000), Barry Johnson (2002), Graham Allan (2007), Ernest Michael
(2013), Charles Read (2015), Jorge Mujica and Marc Thomas (2017),
Berit Stensones and Garth Dales (2022),
Jean Esterle (2023) and Wieslaw \.Zelazko (2025).\\

\noindent {\bf Declaration of interest:} None.\\

\noindent {\bf Competing interest statement:} The author has no
competing interests to disclose.

\end{document}